\documentclass{article}

\usepackage{arxiv}

\usepackage[utf8]{inputenc} 
\usepackage[T1]{fontenc}    
\usepackage{hyperref}       
\usepackage{url}            
\usepackage{booktabs}       
\usepackage{amsfonts}       

\usepackage{amsmath}
\DeclareMathOperator{\arcosh}{arcosh}
\DeclareMathOperator{\arsinh}{arsinh}

\DeclareMathOperator{\gd}{gd}

\usepackage{nicefrac}       
\usepackage{microtype}      
\usepackage{cleveref}       
\usepackage{graphicx}
\usepackage{natbib}
\usepackage{doi}

\title{Catenary and Mercator projection}

\date{\today} 

\author{
    \href{https://orcid.org/0009-0008-6076-6266}{\includegraphics[scale=0.06]{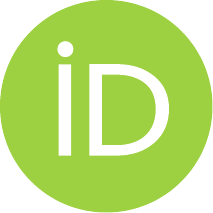}}
    \hspace{1mm} Mikhail A. Akhukov\thanks{Corresponding author: Mikhail Aleksandrovich Akhukov (ma.akhukov@yandex.ru)} \\
     Huawei Nizhny Novgorod Research Center, Nizhny Novgorod, Russia \\
    \texttt{ma.akhukov@yandex.ru} \\
    \And
    \href{https://orcid.org/0000-0002-4930-1846}{\includegraphics[scale=0.06]{orcid.pdf}}
    \hspace{1mm} Vasiliy A. Es’kin \\
     Department of Radiophysics, University of Nizhny Novgorod, Nizhny Novgorod, Russia, 603950 \\
     and \\
     Huawei Nizhny Novgorod Research Center, Nizhny Novgorod, Russia \\
    \texttt{vasiliy.eskin@gmail.com} \\
    \And
    \href{https://orcid.org/0000-0002-0454-5249}{\includegraphics[scale=0.06]{orcid.pdf}}
    \hspace{1mm} Mikhail E. Smorkalov \\
     Skolkovo Institute of Science and Technology, Moscow, Russia \\
     and \\
     Huawei Nizhny Novgorod Research Center, Nizhny Novgorod, Russia \\
    \texttt{smorkalovme@gmail.com}
}

\renewcommand{\shorttitle}{\textit{arXiv} Template}

\hypersetup{
    pdftitle={Catenary and Mercator projection},
    pdfsubject={Differential geometry of hyperbolic functions},
    pdfauthor={Mikhail A. Akhukov, Vasiliy A. Es’kin, Mikhail E. Smorkalov},
    pdfkeywords={Mercator projection, Central cylindrical projection, Gudermannian function, Catenary}
}

\begin{document}
\maketitle

\begin{abstract}

The Mercator projection is sometimes confused with 
another mapping technique, specifically the central cylindrical projection, 
which projects the Earth's surface onto 
a cylinder tangent to the equator, as if a light source is at the Earth's center.
Accidentally, this misconception is rather close to a truth.
The only operation that the map needs is a free bending in 
a uniform gravitational field if the map's material is dense and 
soft enough to produce a catenary profile.
The north and south edges of the map should be parallel and 
placed in the same plane at the appropriate distance.
In this case, the bent map been projected onto this plane gives the Mercator projection.
This property is rather curious, since it allows to make
such a sophisticated one-to-one mapping as the Mercator projection
using simple tools available in the workroom.

\end{abstract}


\section{Introduction}
\label{sec:into}

A catenary is a curve that appears in the profile of an idealized freely hanging heavy chain 
under its own weight, supported only at its ends in a uniform gravitational field. 
The exact shape of the curve is described by the hyperbolic cosine function.

The Mercator projection is a cylindrical map projection. 
It is a conformal projection, which means that it preserves the angles.
The formal equation for the Mercator projection can be expressed in a number of ways.
Within them, there is also a variant that includes a hyperbolic cosine function
\citep{Osborne2016mp}. This fact will be discussed in more details further.

The map in Mercator projection was published in 1569 by 
the Flemish philosopher and mathematician Gerhardus Mercator.
The map was originally designed for navigation, as it maps loxodromes into the straight lines.
The loxodrome or rhumb line is a line on the sphere of a ship's course in the constant compass direction.
This allows sailors to plot a straight course between two points by following a constant compass bearing.
The exact calculation of Mercator projection requires logarithms and modern calculus 
that had not yet been invented at that time.

\begin{figure}[t!]
    \centering
    \includegraphics[width=0.6\textwidth]{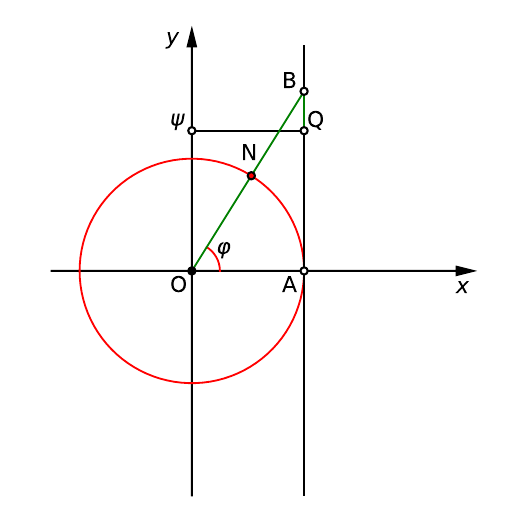}
    \caption{
To calculate the Mercator projection of a given point $N$ with latitude $\varphi$,
we first project the point $N$ to the point $B$ on the vertical line $AB$,
then move it vertically downward to the point $Q$.
The exact vertical coordinate of the point $Q$ is $\psi$,
which is expressed by the equation (\ref{eq:psi_ln}).}
    \label{fig:mercator}
\end{figure}

\section{Mercator projection}
\label{sec:mp}

The Mercator projection is derived from a differential equation,
which is based on the conception of equal scales
and requires the integration of the secant function.
Having a latitude $\varphi$ (see figure \ref{fig:mercator}),
in case of a unit sphere, the formal equation of 
the vertical coordinate of the Mercator map projection is as follow:

\begin{equation}
    \psi := \int\limits_{0}^{\varphi} \frac{dt}{\cos{t}} = 
    \ln \left( \tan \left( \frac{\varphi}{2} + \frac{\pi}{4} \right) \right).
    \label{eq:psi_ln}
\end{equation}

This equation can be expressed in a number of ways using the Gudermannian function
 (\citep{weisstein1999crc} page 778, \citep{Yanpolskiy1960} page 47),
which is popular in calculus, because it relates 
circular ($\varphi$) and hyperbolic ($\psi$) angle measures as $\varphi = gd(\psi)$,
and allows to establish a number of relations between 
trigonometric and hyperbolic function directly without the imaginary unit.

The theory of the Gudermannian function introduces 
a concept of hyperbolic amplitude (\citep{Yanpolskiy1960} page 47),
that is a special angle constructed for an arbitrary point belonging to a unit hyperbola.
It is also known as the gudermannian of the hyperbolic measure $\psi$ because
it is expressed in terms of the Gudermannian function: $\varphi=\gd(\psi)$.

\newpage
\begin{figure}[t!]
    \centering
    \includegraphics[width=0.6\textwidth]{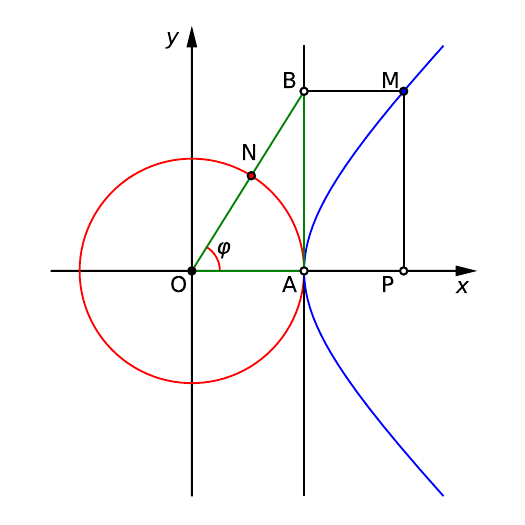}
    \caption{
The hyperbolic amplitude $\varphi$ is constructed for an arbitrary point $M$ 
belonging to the unit hyperbola (blue curve), as shown here.
Also $\varphi$ is the latitude of the point $N$ 
figured out in the Mercator projection (see figure \ref{fig:mercator} for comparison).}
    \label{fig:gudermannian}
\end{figure}

The theory of the Gudermannian function allows to establish a number of relations (see figure \ref{fig:gudermannian}):

\begin{equation}
    L := |AB| = \tan(\varphi) = \sinh(\psi) = |MP|,
    \label{eq:tan_sinh}
\end{equation}

\begin{equation}
    |OB| = \frac{1}{\cos(\varphi)} = \cosh(\psi) = |OP|.
    \label{eq:sec_cosh}
\end{equation}

Using elementary trigonometry in case of $0 \leq \varphi < \pi/2$ we can write:

\begin{equation}
   \tan \left( \frac{\varphi + \pi/2}{2} \right) = 
   \frac{1 + \sin{\varphi}}{\cos{\varphi}} =
   \frac{1}{\cos{\varphi}} + \sqrt{\frac{1}{\cos^2{\varphi}} - 1 }.
   \label{eq:tan_phi_2}
\end{equation}

Equation (\ref{eq:tan_phi_2}) allows to rewrite the right-hand side of equation (\ref{eq:psi_ln}) in the next form:

\begin{equation}
    \psi = 
    \ln \left( \tan \left( \frac{\varphi}{2} + \frac{\pi}{4} \right) \right) =
    \arcosh \left( \frac{1}{\cos(\varphi)} \right).
    \label{eq:psi_ln_arcosh}
\end{equation}

See also equation (\ref{eq:sec_cosh}) for comparison. 
For more details on the integration of the secant function and 
its relation to the history of the Mercator projection see \cite{vis2018history}.

\begin{figure}[t!]
    \centering
    \includegraphics[width=0.6\textwidth]{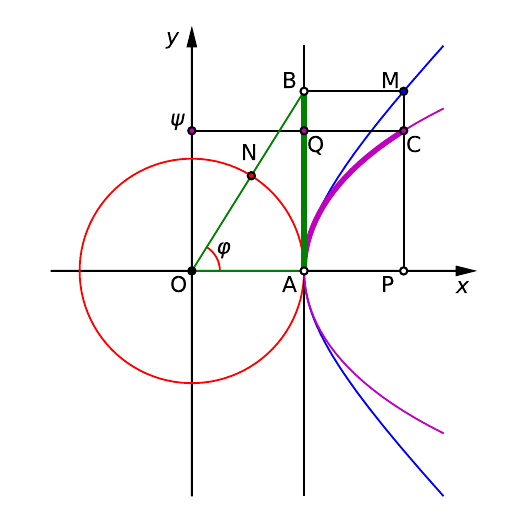}
    \caption{
After introduction of the hyperbolic amplitude $\varphi$ and 
disclose its relation to the hyperbolic measure $\psi$ via the Gudermannian function,
the length of the catenary (solid magenta curve $AC$) introduced implicitly 
via the area hyperbolic cosine as described by equation (\ref{eq:psi_ln_arcosh}) 
is equal to the length $L$ (solid green line $AB$) 
that is vertical component of the central cylindrical map projection. 
}
    \label{fig:catenary}
\end{figure}

Equation (\ref{eq:psi_ln_arcosh}), via area hyperbolic cosine, 
implicitly introduces a catenary oriented in the horizontal direction, 
when vertexes of the unit hyperbola and area hyperbolic cosine touch each other, 
see figure \ref{fig:catenary}.

\newpage

A key point in this work is next.
We need to calculate the length of the area hyperbolic cosine
introduced in such a way as in the equation (\ref{eq:psi_ln_arcosh}),
as a function of the hyperbolic measure $\psi$,
which is the vertical component of the Mercator map projection:

\begin{equation}
    L_{catenary}(\psi) = \\
    \int\limits_{0}^{\psi} \sqrt{1 + \left( \frac{d \cosh(t)} {dt} \right)^2 } dt = \\
    \int\limits_{0}^{\psi} \sqrt{1 + \sinh^2(t) } dt = \\
    \int\limits_{0}^{\psi} \cosh(t) dt = \sinh(\psi).
    \label{eq:L_catenary_psi}
\end{equation}

Using equation (\ref{eq:tan_sinh}) we see that $L_{catenary}(\psi) = \tan(\varphi) = L$,
where $L$ --- is vertical component of the central cylindrical map projection.
The equality of two these lengths, the length of line segment $AB$ and the length of catenary segment $AC$ 
allows to construct the Mercator map projection from the central cylindrical map projection 
by bending the former on the catenary profile and then project the bent map onto the plane.
See figure \ref{fig:catenary}
\footnote{All figures are generated by the scripts in 
Python programming language \citep{VanRossum1995python}
using the matplotlib library \citep{Hunter2007matplotlib}}.

~

~

\subsection{Bending parameter}

Finally, we calculate the exact parameter of the bending for the map in 
the central cylindrical projection to produce the map in the Mercator projection.
If the map in the central cylindrical projection has the plots of a number of loxodromes, 
then the straightening of the loxodromes has to be observed 
from high enough point of view, that will be 
an unavoidable sign of the Mercator map projection.

Suppose, we have a map in the central cylindrical projection constructed between $\alpha$ degrees
of north and south latitudes and physical map's height is $H$. 
The only parameter we need is a distance $D$ between 
the north and south edges of the map been placed parallel and in the same plane.
Using equation (\ref{eq:tan_sinh}) in case of the unit sphere 
the distance $D_{unit~sphere}$ is exactly $2\psi$:
$D_{unit~sphere} = 2\arsinh(\tan(\alpha))$.
In case of sphere with arbitrary radius $R$:

\begin{equation}
    R = \frac{H}{2 \tan{\alpha}},
    \label{eq:R}
\end{equation}
 
we re-scale $D_{unit~sphere}$ by $R$, that gives the value of distance $D$:

\begin{equation}
    D = \frac{H}{\tan{\alpha}} \arsinh \left( \tan{\alpha} \right).
    \label{eq:D}
\end{equation}

\section{Conclusion}

We show that the map in the Mercator projection can be 
constructed from the map in the central cylindrical projection
as if some one has such map been printed on dense and soft enough material
allows it to hang freely to make the catenary profile and then project the bent map onto the plane.
Interestingly speaking, if such map has a number of loxodromes, then 
a person standing high enough can see that the loxodromes become 
straight lines when original map is bent when, for example, 
the map is carried by two other persons if they hold the map at the north and south edges.

This hypothetical situation could happen in real life 
in some kind of workroom and serve as 
a basis for further study of this property ---
the straightening of loxodromes printed on the map in 
the central cylindrical projection when this map hangs freely as described above.

\bibliographystyle{unsrtnat}
\bibliography{references}

\end{document}